\documentclass[12pt]{amsart}
\usepackage{amssymb}
\usepackage{relsize}

\newcommand{\st}{_{st}}

\newcommand{\cho}{\mtc{C}o}
\newcommand{\mtjccp}{\mtc J_{\ccp}}

\newcommand{\mtjce}{\mtj(\ce)}
\newcommand{{\bm}}{{\overline{m}}}

\newcommand{\almcd}{\lamcd}

\newcommand{\sccj}{{|\ccj|}}
\newcommand{\dimccj}{\dim(\cc^j)}

\newcommand{\ccp}{{\cc'}}
\newcommand{\gr}{\mtr{Gr}}

\newcommand\numeq[1]%
  {\stackrel{\scriptscriptstyle(\mkern-1.5mu#1\mkern-1.5mu)}{=}}

\newcounter{relctr} 
\everydisplay\expandafter{\the\everydisplay\setcounter{relctr}{0}} 

\AtBeginDocument{} 

\newcommand{\barcj}{\mtr C_j}

\newcommand{\io}{{i_1}}
\newcommand{\itw}{{i_2}}

 \newcommand{\phibr}{\phi_{\bar R}}
 \newcommand{\jc}{{\mtc J}^{c}}

\usepackage{tikz-cd}
\usepackage{latexsym,amssymb,amsmath}
\usepackage{hyperref}
\usepackage{comment}
\usepackage{mathrsfs}
\usepackage{amsmath,amscd}
\usepackage[enableskew]{youngtab}
\usepackage[all]{xy}
\usepackage{empheq}

\usepackage{amsmath}

\usepackage{empheq}
\usepackage{xcolor}
\definecolor{lightgreen}{HTML}{90EE90}

\newcommand{\irrcc}{\irr(\cc)}
\newcommand{\irrcd}{\irr(\cd)}

\newcommand{\mtj}{{\mtr{supp}}}

\newcommand{\lag}{\langle}
\newcommand{\rag}{\rangle}

 \newcommand{\vsk}{\vskip 0.15cm \noindent}

\newcommand{\ccj}{\cc_j}

 \newcommand{\fpcc}{\fp(\cc)}

\newcommand{\vect}{\mtr{Vec}}

\newcommand{\sent}{\mapsto}

\newcommand{\muj}{\mu_j}

\numberwithin{equation}{section}

\newcommand\C{\mathcal{C}}
\DeclareMathOperator{\id}{id}

\excludecomment{verlong}
\includecomment{vershort}
\excludecomment{noncompile}

\usepackage{tikz}
\usetikzlibrary{matrix}
\usepackage{combelow}
\newcommand{\ccb}{{\mathcal B}}

\usepackage{amsmath,amsthm,amsfonts,amssymb}
\usepackage{amsmath, amsthm, amssymb, amscd, amsfonts}
\usepackage[all]{xy}
\usepackage{amssymb, amsthm, amsmath}
\usepackage{amsfonts}
\usepackage{color}

\newcommand{\ra}{\rightarrow}

\newcommand{\ot}{\otimes}
\newcommand{\co}{\mathcal O}
\newcommand{\xra}{\xrightarrow}
\newcommand{\mtc}{\mathcal}
\newcommand{\cs}{\mtc S}

\newcommand{\lam}{\lambda}

\newcommand{\Lam}{\Lambda}

\newcommand{\al}{\alpha}
\newcommand{\eps}{\epsilon}

\newcommand{\ul}{\underline}

\numberwithin{equation}{section}
\newtheorem{lem}[equation]{Lemma}

\theoremstyle{plain}

\newtheorem{thm}[equation]{Theorem}
\newtheorem{prop}[equation]{Proposition}
\newtheorem{defn}[equation]{Definition}
\newtheorem{cor}[equation]{Corollary}
\newtheorem{rem}[equation]{Remark}
\newcommand{\ch}{\chi}
\newcommand{\mtr}{\mathrm}

\numberwithin{equation}{section}
\newcommand{\ncm}{\newcommand}
\ncm{\np}{\newpage}
\ncm{\ebl}{\end{thebibliography}}
\ncm{\bbl}{\begin{thebibliography}}
\ncm{\chd}{_{ _{\ch}}}
\ncm{\ald}{_{ _{\al}}}
\newcommand{\blam}{\Lam}
\ncm{\cP}{\mathcal{P}}
\ncm{\ei}{e_i}
\ncm{\eij}{e_{i,\;j}}
\ncm{\bt}{\begin{thm}}
\ncm{\bdef}{\begin{defn}}
\ncm{\edf}{\end{defn}}
\ncm{\et}{\end{thm}}
\ncm{\bc}{\begin{cor}}
\ncm{\bl}{\begin{lem}}
\ncm{\el}{\end{lem}}
\ncm{\bpf}{\begin{proof}}
\ncm{\epf}{\end{proof}}
\ncm{\ec}{\end{cor}}
\ncm{\ord}{\mtr{ord}}
\ncm{\er}{\end{rem}}
\ncm{\br}{\begin{rem}}
\ncm{\bn}{\begin}

\ncm{\bp}{\begin{prop}}
\ncm{\ep}{\end{prop}}
\ncm{\bd}{
\begin{document}}
\ncm{\ed}{\end{document}}
\ncm{\beq}{\begin{equation}}
\ncm{\beqn}{\begin{equation*}}
\ncm{\eeq}{\end{equation}}
\ncm{\eeqn}{\end{equation*}}
\ncm{\bea}{\begin{eqnarray}}
\ncm{\eea}{\end{eqnarray}}
\ncm{\beanon}{\begin{eqnarray*}}
\ncm{\eeanon}{\end{eqnarray*}}\ncm{\ek}{\eps|_K}\ncm{\diez}{\#}
\ncm{\bwt}{\bowtie}
\ncm{\cC}{\mtc{C}}\ncm{\cc}{\mtc{C}}
\ncm{\cX}{\mtc{X}}
\ncm{\wt}{\widetilde}
\ncm{\sg}{\sigma}
\ncm{\Rep}{\mathrm{Rep}}
\DeclareMathOperator{\Irr}{Irr}
\ncm{\X}{\mathcal{X}}
\ncm{\cA}{\mathcal{A}}
\ncm{\HKer}{\mtr{HKer}}
\ncm{\LKER}{\mtr{LKer}}
\ncm{\aad}{\mtr{ad}}
\newcommand{\mbf}{\mathbb F}
\ncm{\Dr}{\mtr{D}}
\ncm{\cD}{{\mathcal{D}}}\ncm{\cd}{{\mathcal{D}}}\ncm{\ce}{{\mathcal{E}}}
\ncm{\G}{\mathcal{G}}
\ncm{\Dc}{\mtc{D}}
\ncm{\E}{\mtc{E}}
\ncm{\fp}{\mtr{FPdim}}
\ncm{\Vc}{\mtr{Vec}}
\ncm{\cK}{\mtc{K}}
\ncm{\cM}{\mtc{M}}
\ncm{\cE}{\mtc{E}}
\ncm{\cS}{\mtc{S}}

\newcommand{{\ipr}}{i'}

\DeclareMathOperator{\End}{End}
\ncm{\cop}{\mtr{cop}}
\ncm{\op}{\mtr{op}}
\ncm{\chr}{character }\ncm{\ck}{\mtc{K}}
\ncm{\bw}{\bwt}
\ncm{\hker}{\mtr{HKer}}
\ncm{\bx}{\boxtimes}
\ncm{\blue}{\textcolor[rgb]{.00, .00, 1.00}}
\ncm{\red}{\textcolor[rgb]{1.00, .00, .00}}
\ncm{\green}{\textcolor[rgb]{.50, 0.20, .90}}
\ncm{\bne}{\begin{enumerate}}
\ncm{\ene}{\end{enumerate}}
\ncm{\lker}{\mtr{LKer}}
\ncm{\md}{\medbreak}
\ncm{\rep}{\Rep}\ncm{\ind}{\mtr{ind}}
\ncm{\mdn}{\md\noindent}
\ncm{\dd}{$}
\ncm{\up}{^}
\newcommand{\tcs}{\text}
\newcommand{\mbb}{\mathbb B}
\newcommand{\vs}{\mathbb V}
\newcommand{\sth}{suppose that\;}
\newcommand\rad{\operatorname{rad}}
\newcommand{\itm}{\item}
\newcommand{\dbd}{$$}
\newcommand{\mol}{\mtr{mod}}
 \newcommand{\ro}{\rho}
\newcommand{\irr}{\mathrm{Irr}}
\newcommand{\mbc}{\mathbb C}
\newcommand{\mbs}{\mathbb S}
\newcommand{\mbz}{\mathbb Z}
\newcommand{\ct}{\mtc T}
\newcommand{\sm}{\setminus}
\newcommand{\epl}{^{+}}
\newcommand{\sbsq}{\subseteq}
\newcommand{\sbs}{\subset}
\newcommand{\cco}{\mtr{co}}
\newcommand{\cz}{\mathcal{Z}}
\newcommand{\dual}{^{*}}
\newcommand{\Gm}{\Gamma}
\ncm{\cY}{\mtc{Y}}
\newcommand\ZZ{{\mathbb Z}} 
\newcommand{\bab}{\color{DarkOrchid}{}}
\newcommand{\eab}{\normalcolor{}}
\newcommand{\subs}{\subsection}
\newcommand{\cv}{\mtc{V}}
  \newcommand{\grn}{\green}
\newcommand{\dt}{\delta}

\newcommand{\ccf}{\mathrm{ {CF}(\cc)}}
\newcommand{\cce}{\mathrm{ {CE}(\cc)}}
\newcommand{\cecc}{\mathrm{ {CE}(\cc)}}
\newcommand{\cecd}{\mathrm{ {CE}(\cd)}}
\newcommand{\kk}{\Bbbk}
\newcommand{\otL}{\ot_{L}}
\newcommand{\otl}{\ot_{L}}
\newcommand{\unpsi}{1_{\psi}}
\newcommand{\epsi}{e_{\psi}}
\newcommand{\ephi}{e_{\phi}}
\newcommand{\ech}{e_{\ch}}
\newcommand{\nleftcid}{\text{left normal  coideal subalgebra}}
\newcommand{\dimL}{\dim_{\kk}L}
\newcommand{\cl}{\mtc L}
\newcommand{\mj}{\mtc J}
\newcommand{\tl}{\tilde L}
\newcommand{\tL}{\tilde L}
\newcommand{\tpsi}{\tilde(\psi)}
\newcommand{\tmx}{\tilde{\mtc X}}
\newcommand{\zlh}{\mathrm{ZL}}
\newcommand{\ba}{\mathrm A}
\newcommand{\bv}{\mathrm V}
\newcommand{\zhopf}{\mtc{Z}_{\mtr{Hopf}}}
\newcommand{\lstar}{L^{*}}
\newcommand{\ldstar}{L^{**}}
\newcommand{\mstar}{M^{*}}
\newcommand{\mdstar}{M^{**}}
\newcommand{\lkera}{\lker_{A}}
\newcommand{\mdprime}{M''}
\newcommand{\ldprime}{L''}
\newcommand{\cm}{\mtc M}
\newcommand{\ccm}{\mathcal M}
\newcommand{\cn}{\mathcal N}
\newcommand{\ccn}{\mathcal N}
\newcommand{\rx}{\mtr{Rex}}
\newcommand{\cca}{\ca}
\newcommand{\ih}{\underline{\mtr{Hom}}}
\newcommand{\cih}{\underline{\mtr{coHom}}}
\newcommand{\hm}{\mtr{ {Hom}}}
\newcommand{\cov}{\mtr{coev}}
\newcommand{\rora}{\rho^{\mtr{ra}}}
\newcommand{\rola}{\rho^{\mtr{la}}}
\newcommand{\cx}{\mtc X}
 \newcommand{\cZ}{\cz}
 \newcommand{\ca}{\cA}
 \newcommand{\stat}{\noindent}
 \newcommand{\bfa}{{\bf A}}
 \newcommand{\unu}{\mathbf{1}}
 \newcommand{\barzu}{{\bar {  Z}(\unu)}}
 
\newcommand{\idx}{\id_X}
\newcommand{\lprime}{L'}
\newcommand{\mprime}{M'}
\newcommand{\nat}{ \mtr{{  Nat}}}
\newcommand{\ft}{\mtc F_\lam}
\newcommand{\rhau}{\rightharpoonup}
\newcommand{\lhau}{\leftharpoonup}
\newcommand{\cf}{\mathrm{ {CF}}}

\newcommand{\cfc}{\mathrm{{CF}}(\cc)}
\newcommand{\csu}{\overline{\mathfrak{  C}}}
\newcommand{\cfcc}{\mathrm{ {CF}}(\cc)}
\newcommand{\cfcd}{\mathrm{CF}(\cd)}
\newcommand{\cfd}{\mathrm{CF}(\cd)}
\newcommand{\czcc}{{\cz(\cc)}}
\newcommand{\czcd}{{\cz(\cd)}}
\newcommand{\czt}{{\cz(\cz(\cc))}}
\newcommand{\enx}{\mtr{  End}}
\newcommand{\runu}{R(\unu)}

\newcommand{\bdfn}{\bn{defn}}
\newcommand{\edfn}{\end{defn}}
\newcommand{\deltax}{\delta_X}
\newcommand{\deltav}{\delta_V}
\newcommand{\repcca}{\rep_\cc(A)}
\newcommand{\xotay}{X \ot_A Y}
\newcommand{\xoty}{X \ot Y}
\newcommand{\votw}{V \ot W}
\newcommand{\votaw}{V \ot_A W}
\newcommand{\dimax}{\dim_AX}
\newcommand{\dimccx}{\dim_\cc(X)}
\newcommand{\dimcca}{\dim_\cc(A)}
\newcommand{\dimccv}{\dim_\cc(V)}
\newcommand{\dima}{\dim_A}
\newcommand{\biga}{A}
\newcommand{\comp}{\mathbb C}
\newcommand{\tehtaa}{\theta_A}
\newcommand{\tetaa}{\theta_A}
\newcommand{\ida}{\id_A}
\newcommand{\hma}{\hm_A}
\newcommand{\hmcc}{\hm_\cc}
\newcommand{\fv}{F(V)}
\newcommand{\fw}{F(W)}
\newcommand{\ota}{\ot_A}
\newcommand{\repza}{\rep_\cc^0(A)}
\newcommand{\epsa}{\eps_A}
\newcommand{\bndefn}{\bn{defn}}
\newcommand{\edefn}{\end{defn}}
\newcommand{\bdefn}{\bn{defn}}

\newcommand{\vld}{V^{*}}
\newcommand{\vldd}{V^{**}}
\newcommand{\xld}{X^{*}}
\newcommand{\xldd}{X^{**}}
\newcommand{\yld}{Y^{*}}
\newcommand{\yldd}{Y^{**}}
\newcommand{\aldu}{A^{*}}
\newcommand{\aldd}{A^{**}}

\newcommand{\ia}{\mtr{i}_A}
\newcommand{\aota}{A\ot A}

\newcommand{\idv}{\id_V}

\newcommand{\ld}{^*}
\newcommand{\repg}{\rep(G)}

\newcommand{\thetav}{\theta_V}

\newcommand{\tta}{\theta_A}

\newcommand{\muv}{\mu_V}
\newcommand{\muw}{\mu_W}

\newcommand{\dimcc}{\dim(\cc)}
\newcommand{\chii}{\chi_i}
\newcommand{\chistar}{\ch_{i^*}}
\newcommand{\chj}{\ch_j}
\newcommand{\chm}{\ch_m}
\newcommand{\chn}{\ch_n}
\newcommand{\dimvi}{\dim(V_i)}
\newcommand{\mtcd}{Q}
\newcommand{\mtca}{\mtc A}
\newcommand{\lamcd}{\lam_\cd}
\newcommand{\fpdimcd}{\fp(\cd)}
\newcommand{\laml}{\lam_L}
\newcommand{\apm}{A//M}
\newcommand{\apl}{A//L}
\newcommand{\repapm}{\rep(\apm)}
\newcommand{\repapl}{\rep(\apl)}
\newcommand{\dimvj}{\dim(V_j)}
\newcommand{\dvi}{\dim(V_i)}
\newcommand{\dvj}{\dim(V_j)}
\newcommand{\sumjtom}{\sum_{j=0}^m}
\newcommand{\sumitom}{\sum_{i=0}^m}
\newcommand{\sij}{s_{ij}}
\newcommand{\sji}{s_{ji}}
\newcommand{\dxj}{d_j}
\newcommand{\dxi}{d_i}
\newcommand{\dimka}{\dim_{\kk}(A)}
\newcommand{\dimk}{\dim_{\kk}}
\newcommand{\blaml}{\blam_L}
\newcommand{\sumjtor}{\sum_{j=0}^r}
\newcommand{\dimkl}{\dim_{\kk}(L)}
\newcommand{\mtcjl}{\mtc J_L}
\newcommand{\vota}{ V\ot A}
\newcommand{\vi}{V_i}
\newcommand{\vj}{V_j}
\newcommand{\dimcd}{\dim(\cd)}

\newcommand{\alij}{\al_{ij}}
\newcommand{\alji}{\al_{ji}}
\newcommand{\rcc}{r_\cc}
\newcommand{\rcd}{r_\cd}
\newcommand{\clsx}{[X]}
\newcommand{\clsy}{[Y]}
\newcommand{\clsz}{[Z]}
\newcommand{\rcdp}{r_{\cd'}}
\newcommand{\sumjtorp}{\sum_{j=0}^{r'}}
\newcommand{\aljm}{\al_{jm}}
\newcommand{\aljn}{\al_{jn}}
\newcommand{\sjm}{s_{jm}}
\newcommand{\smj}{s_{mj}}
\newcommand{\snj}{s_{nj}}

\newcommand{\betaij}{\beta_{ij}}
\newcommand{\betaji}{\beta_{ji}}

 \newcommand{\ip}{i'}
\newcommand{\sumjtoprp}{\sum_{j=0}^{r'}}
\newcommand{\sumjtopr}{\sum_{j=0}^{r}}
 \newcommand{\teh}{\tilde{h}}
\newcommand{\cdp}{{\cd'}}\newcommand{\xphii}{X_{\phi(i)}}
\newcommand{\inv}{^{-1}}

\newcommand{\fq}{f_Q}
\newcommand{\tr}{\mtr{tr}}
\newcommand{\rtwone}{R_{21}R}

\newcommand{\ccad}{{\cc_{\mtr{ad}}}}
\newcommand{\ccpt}{{\cc_{\mtr{pt}}}}
\newcommand{\qtr}{quasi-triangular\;}
\newcommand{\trq}{\tr_q}

\newcommand{\repal}{\mtr{Rep}(A//L)}
\newcommand{\lkeravi}{\lker_A(V_i)}
\newcommand{\lkeravj}{\lker_A(V_j)}
\newcommand{\cross}[1][1pt]{\ooalign{%
 \rule[1ex]{1ex}{#1}\cr
 \hss\rule{#1}{.7em}\hss\cr}}
\newcommand{\blml}{\blam_L} 
\newcommand{\phir}{\phi_R}
\newcommand{\kda}{{  \Phi(A)}}

\newcommand{\mtcil}{\mtc{I}_L}

\newcommand{\un}{\unu}
\newcommand{\tfl}{\mtc{T}}
\newcommand{\barzm}{\barz(M)}
\newcommand{\barzn}{\barz(N)}
\newcommand{\ccr}{\mtc R^{\cc}}
\newcommand{\ulc}{\ul{\cc}}

\newcommand{\pimx}{\pi_{M;\;X}}
\newcommand{\pinx}{\pi_{N;\;X}}
\newcommand{\acc}{{\mathrm A_\cc}}
\newcommand{\epsu}{\eps_\unu}

\newcommand{\ob}{\mtr{Obj}}
\newcommand{\obc}{\mtr{Obj(\cc)}}
\newcommand{\ccop}{\cc^{\mtr{op}}}
\newcommand{\mtf}{\mtc F_\lam}
\newcommand{\mtfi}{\mtc F^{-1}_\lam}
\newcommand{\elcd}{\ell_\cd}
\newcommand{\mcid}{\mtc I_\cd}
\newcommand{\mcidp}{\mtc I_{\cd'}}
\newcommand{\wtildelcd}{\widetilde{\elcd}}
\newcommand{\wtildelcdp}{\widetilde{\ell_{\cd'}}}
\newcommand{\cpt}{\cc_{\mtr{pt}}}
\newcommand{\barzr}{\barz_\cd}
\newcommand{\barzv}{\barz(V)}
\newcommand{\acd}{\mathrm A_\cd}
\newcommand{\czrcd}{\cz_\cc(\cd)}
\newcommand{\sml}{\Small}
\newcommand{\bs}{\blue{\Small }}
\newcommand{\yd}{Yetter-Drinfeld}

\newcommand{\sumitor}{\sum_{i=0}^r}
\newcommand{\cdop}{\cd^{\mtr{op}}}
\newcommand{\ccrev}{\cc^{\mtr{rev}}}
\newcommand{\barz}{{\bar{\mathrm Z}}}
\newcommand{\etl}{etale\;}
\newcommand{\czca}{\cz(\ca)}

\newcommand{\lkerccv}{{\lker_\cc(V)}}
\newcommand{\ml}{{\mtr L}}
\newcommand{\ma}{{\mtr A}}
\newcommand{\barm}{{\frac{R_m}{\fp(R_m)}}}
\newcommand{\mk}{{\mtr K}}
\newcommand{\mce}{\mtr{CE}}
\newcommand{\mh}{\mtc H}
\newcommand{\cccd}{{\big(\cc\big/ \cd\big)_r}}

\bd
\newcommand{\ccpp}{{\cz_2(\cc)}}

\title[Braided fusion categories]{On some Frobenius type divisibility results in a  premodular category}
\author{Sebastian Burciu}
\address{Inst.\ of Math.\ ``Simion Stoilow" of the Romanian Academy P.O. Box 1-764, RO-014700, Bucharest, Romania}
\email{sebastian.burciu@imar.ro}
\date{\today}
\maketitle
\section*{Abstract}
In this note some new Frobenius type divisibility results are obtained for premodular categories. In particular,  we extend \cite[Corollary 3.4]{yu-sd} from the settings of  super-modular categories  to arbitrary pseudo-unitary premodular categories.

\section{Introduction}
In \cite{yu-sd} the author showed that in a slightly-degenerate braided fusion category one has that $\frac{\dimcc}{2\dim(Y)^2}\in \mathbb A$ for any simple object $Y$ of $\cc$. The main goal of this note is to generalize this result in the settings of a premodular fusion category with pointed M\"uger center, i.e. all the objects of this  center are invertible.

Our first main result is the following: 

\bt\label{squaredim} 
Let $\cc$ be a premodular category and $\cd$ a fusion subcategory of $\cc$. If $\cd\cap\ccpp=\vect$ then 
\beq\label{8}
\frac{\dimcc}{\dim(Y)^2}\in \mathbb A
\eeq
for any simple object $Y$ in $\cd$.
\et

Moreover,  in the case of a pointed M\"uger center we prove that:
 \bt\label{main3}
Let  $\cc$ be a pseudo-unitary premodular  fusion category with $\cz_2(\cc)\subseteq \ccpt$.
\bne
\item
For any simple object $\cc$ of $Y$ one has
\beq\label{ptdiv}
\frac{\fp(\cc)\fp(\cz_2(\cc))}{\fp(Y)^2}\in \mathbb A.
\eeq
\item 
If the action of $\cz_2(\cc)$ on the simple objects of $\cc$ is free then 
\beq\label{freeptdiv}
\frac{\fp(\cc)}{\fp(\cz_2(\cc))\fp(Y)^2}\in \mathbb A.
\eeq
\ene
\et
Recall that $\ccpt$ denotes the maximal pointed subcategory of $\cc$. Note that the second item of the above Theorem generalizes \cite[Corollary 3.4]{yu-sd} in the case of a pseudo-unitary premodular category.

Shortly, this note is organized as follows. In Section \ref{prelim} we recall the basic properties of fusion categories needed through this paper. In Subsection \ref{adj} we briefly recall the adjoint algebra and the characters of simple objects of pivotal fusion categories as defined in \cite{scalg}. 

In Section \ref{hecke} we construct a Hecke type algebras associated to the right cosets of a fusion subcategory of a given fusion category. In the same section we prove orthogonality relations for Hecke algebras. In Section \ref{mc} we prove the main results stated above.
\section{Preliminaries}\label{prelim}
In this section we recall the main properties of premodular fusion categories that are needed through the paper.  For the basic theory of fusion categories, we refer the reader to \cite{EGNO15}. Recall that a fusion category is by definition, a semisimple finite tensor category.

Throughout this note $\cc$ denotes a fusion category over $\comp$ with eventual additional properties. Let $\irr(\cc):=\{X_0, X_1, \dots, X_m\}$ be a complete set of representatives for the isomorphism classes of simple objects of $\cc$. By $K_0(\cc)$ we denote the Grothendieck ring of $\cc$  and let $K(\cc):=K_0(\cc)\ot_{\mathbb Z}\comp$.

Recall that a {\it pivotal structure} of a rigid monoidal category $\cc$ is an isomorphism $j:\id_\cc\ra ()^{**}$ of monoidal functors. A pivotal monoidal category is a rigid monoidal category endowed with a pivotal structure. For any pivotal structure one can construct categorical dimensions $\dim(X)$ of any object $X$ of $\cc$. In fact, the map $\dim:K(\cc)\ra\comp, [X]\sent \dim(X)$ is an algebra morphism.

A pivotal structure a on a fusion category $\cc$ is called {\it spherical}, see \cite{EGNO15} if
$\dim(V)=\dim(V^*)$ for any object $V$ of $\cc$. A tensor category is spherical if it is
equipped with a spherical structure.
Given a fusion category there is a unique algebra morphism $\fp:K(\cc)\ra\comp$ such that $\fp(X_i)>0$ for any $X_i\in \irr(\cc)$. Then $\fp$ is called the Frobenius-Perron morphism.

Recall that a braided fusion category is called {\it premodular} if it has a spherical structure. Equivalently, a premodular category  is a braided  fusion category equipped with a twist
(also called a balanced structure), see \cite[Section 8.10 ]{EGNO15} for details.

A fusion category $\cc$ is called {\it pseudo-unitary} if $\fpcc =\dimcc$. If such is the case,
then by \cite[Proposition 8.23]{ENO}, $\cc$ admits a unique spherical structure with respect to which the categorical dimensions of  each simple coincides with its Frobenius-Perron dimension, i.e. $\fp(X)=\dim(X)$ for any object $X$ of $\cc$. 
\subsection{Adjoint algebra and internal characters}\label{adj}
Let $\cc$ be a fusion category and $\czcc$ be its Drinfeld center. The forgetful functor $F:\czcc\ra\cc$  admits a right adjoint functor $R:\cc \ra \czcc$  and let  $Z :=FR:\cc \ra \cc$. Then $Z$ is a Hopf comonad and by \cite[Section 2.6]{scalg} one has that 
\beq
Z(V)\simeq \int_{X\in \cc}X\ot V\ot X^*
\eeq
It is known that $A:=R(\unu)$ has the structure of central commutative algebra in $\cz(\cc)$.

We denote by $m:A\ot A \ra A$ the multiplication  and by $u:\unu\ra A$ the unit of the
adjoint algebra $A = Z (\unu)$.

The vector space $\cecc:= \hm_{\C}(\unu, A) $ is called {\it the set of central elements.} For $a, b \in \cecc$ , we set $a ·b :=m \circ (a \ot b)$. Then the set $\cecc$ is a $\comp$-algebra with respect to this operation. 
The space $\cfcc:=\hm_\cc(A, \unu)$ is called the {\it space of class functions} of $\cc$.  The space $\cfcc$ is a $\comp$-algebra with the multiplication defined $f\star g:=f \circ Z(g) \circ \delta_{\unu}.$ for any  $f, g\in \cfcc$.   Here $\delta: Z \ra  Z^2$ is the comultiplication structure of the Hopf comonad $Z$, see \cite{scalg} for details.  To any object $X $ of $\cc$ Shimizu assigned in \cite[Section 3]{scalg} a class function $\mtr{ch}(X)$ called {\it the internal character}. It is defined as a partial trace of the canonical left action $\ro_X:A\ot X\ra \unu$.  By \cite[Theorem 3.10]{scalg} one has that $\mtr{ch}(X\ot Y)=\mtr{ch}(X)\mtr{ch}(Y)$ for any two objects $X$ and $Y$ of $\cc$ and $K(\cc)\simeq \cfcc$ as algebras. 
 
There is a also non-degenerate a paring
$\langle\;,\; \rangle : \cfcc \times \cecc\ra \unu$, given by $ \langle f, a\rangle  \id_{\unu}= f \circ a,$ for all $f \in \cfcc$ and $a\in \cecc$. 

For any simple object $X_i$ of $\cc$ we denote by $\ch_i:=\mtr{ch}(X_i)$ its associated character. The central element space $\cecc$ has a basis of primitive orthogonal idempotents $E_i$ such that   $<\ch_i, E_j>=\delta_{i,j}d_i$ for all $i,j$.

Recall $R:\cc \ra \czcc$ is a right adjoint to the forgetful functor $F:\czcc \ra \cc$. As explained in \cite[Theorem 3.8]{scalg} this adjunction  gives an isomorphism of algebras
\beq\label{adjisom}
\cfcc \xra{\cong} \mtr{End}_{\czcc}(R(\unu)),\;\; \ch\mapsto Z (\ch)\circ \delta_\unu.
\eeq
Since $\czcc$ is also  fusion category we can write $R(\unu)=\bigoplus_{j=0}^m\mathcal C^j$ as a direct sum of simple objects in $\czcc$. Recall that $\cc^j$ are called {\it conjugacy classes } for $\cc$.

For the rest of this section let $\cc$ be a fusion category with a commutative Grothendieck ring $K(\cc)$. Since $\cfcc$ is also commutative semisimple $\comp$-algebra has
\beq\label{cfdec}
\cfcc\simeq \bigoplus_{j=0}^m \comp F_j
\eeq
where $F_j$ are the central primitive idempotents of $\cc$. We define $\mtc J:=\{0,\dots, m\}$.

Since $\cfcc$ is a commutative algebra it follows from Equation \eqref{adjisom} that $R(\unu)$ is multiplicity free, i.e. each conjugacy class $\cc^j$ appears with multiplicity $1$ inside $R(\unu)$. Moreover we fix the canonical  bijection $F_j \leftrightarrow\cc^j$ for which the image of $F_j$ under the canonical adjunction isomorphism from Equation \eqref{adjisom} is  the projection of $R(\unu)$ on $\cc^j$.

With this bijection, the second orthogonality from \cite[Theorem 6.10]{scalg} can be written as:
\beq\label{somu}
\sum_{ i=0}^m\mu_l(\ch_i)\mu_k(\ch_{i^*})=\delta_{l,k}\frac{\dimcc}{\dim(\cc^k)}
\eeq
where $\muj:\cfcc\ra\comp$ is the set of algebra maps of $\cfcc$. In terms of Equation \eqref{cfdec} one has $\muj(F_k)=\delta_{j,k}$.

Following \cite{ccc-march} we use the following definition:
\bn{defn} To any fusion subcategory $\cd$ we associate a set of indices $\mtc J_\cd\subseteq \mtc J$ such that:
\beq\label{lamcsl}
\lam_{\cd}=\sum_{j\in \mtc J_\cd}F_j.
\eeq
\end{defn}

In this settings, Equation (4.20) from \cite{ccc-march} gives that
\beq\label{sumdim}
\sum\limits_{j \in \mtc J_{\cd}}\dim(\cc^j)=\frac{\dim(\cc)}{\dim(\cd)}.
\eeq
\section{Hecke algebras of cosets}\label{hecke}
For  a fusion subcategory $\cd$ of $\cc$ we denote by $\simeq^r_\cd$ the equivalence relation on the set of simple objects $\irrcc$ of $\cc$. It is given by $X\simeq^r_\cd Y$ if and only if there is a simple object $S$ of $\cd$ such that $X$ is a constituent of $Y\ot S$. An equivalence class of $\simeq^r_\cd$ is called a {\it right coset} of $\cc$ with respect to $\cd$.  We denote by $\cccd$ the set of equivalence classes (right cosets) with respect to $\cd$.  It was shown in \cite{bbr} that 
\beq\label{initcoset}
X\simeq^r_{\cd}Y\iff \frac{[X]R_\cd}{\fp(X)}= \frac{[Y]R_\cd}{\fp(Y)}
\eeq
where $R_\cd$ is the regular element of the fusion subcategory $\cd$, $R_\cd=\sum_{X\in \irrcd}\fp(X)[X]\in K(\cc)$. In this case it can be shown that
\beq\label{cformula}
\frac{[X]R_\cd}{\fp(X)}= \frac{[Y]R_\cd}{\fp(Y)}=\frac{R_m}{\fp(R_m)}
\eeq
where $m\in \cccd$ is the right coset corresponding to $X$ and $Y$ and and $R_m:=\sum_{Z\in m}\fp(Z)[Z]\in K(\cc)$ is the regular element corresponding to any coset $m$.

Based on Equation \eqref{cformula} it is easy to verify that if $K(\cc)$ is commutative then for all $m,n,p\in \cccd$
\beq
\frac{R_m}{\fp(R_m)}\frac{R_n}{\fp(R_n)}=\sum_{p\in \cccd} H^p_{mn} \frac{R_p}{\fp(R_p)}
\eeq
where $H^p_{mn}$ are (integral) scalars given by $$H^p_{mn}=\big(\sum_{Z\in p}\fp(Z)N^Z_{XY}\big)$$
where  $X\in m$ and $Y\in n$ are arbitrarily chosen representatives for each coset $m, n \in \cccd$. Thus in the case $K(\cc)$ is commutative, the vector subspace $\mh(\cd):=\{R_m\;\big|\; m\in \cccd\}$ of $K(\cc)$ generated by the regular elements $R_m$ is in fact a subalgebra of $K(\cc)$. It is called the {\it Hecke algebra} of $\cc$ with respect to $\cd$. Note that Equation \eqref{cformula} implies $\mh(\cd)=K(\cc) R_{\cd}$ inside $K(\cc)$. 

 \bp\label{dimhecke} 
 Let $\cd$ be a fusion subcategory of a pseudo-unitary fusion category $\cc$ with a commutative Grothendieck ring.  Then the dimension of the Hecke algebra  $\mh(\cd)$ equals the cardinality of $\mtc J_\cd$. 
 \ep
 \bpf
 If $\cc$ is pseudo-unitary then as above $K(\cc)\simeq \cfcc$ and $R_\cd$ corresponds to the idempotent cointegral $\lamcd$. Then from Equation \eqref{cfdec} in this case one has:
 $$\mh(\cd)=\cfcc\lam_\cd=\bigoplus_{j \in \mtc J_\cd} \cfcc F_j=\bigoplus_{j \in \mtc J_\cd} \comp F_j
 $$
 by Equation \eqref{lamcsl}.
 \epf
  \subsection*{General Hecke orthogonality}

For any $t\in \cccd$ let $t^*:=\{X^*\;|\;X\in t\}.$ It is easy to see that in the case of a commutative Grothendieck ring, if  $t\in \cccd$ is a right coset then $t^*$ is also a right coset. Indeed, since $K(\cc)$ is commutative one has $([X])R_\cd)^*=R_\cd^*[X]^*=R_\cd[X^*]$ since $R_\cd^*=R_\cd$.
\newcommand{\xit}{X_t}
\newcommand{\xitstar}{{X_{t^*}}}
\newcommand{\xis}{X_s}
\newcommand{\xisstar}{{X_{s^*}}}
\newcommand{\chit}{\ch_t}
\newcommand{\chitstar}{\ch_{t^*}}
\newcommand{\chisstar}{\ch_{s^*}}
\newcommand{\chis}{{\ch_{ _{i_s}}}}
\newcommand{\hacd}{\mh(\cd)}
For any $t\in \cccd$ choose a representative simple object $\xit$ of $\cc$ belonging to this coset such that $\xit^*=\xitstar$. Denote also $\chit:=\ch(\xit)\in \cfcc$.
 
Equation \eqref{lamcsl} implies that $$\hacd=\cfcc\almcd=\bigoplus_{j\in \mtc J_\cd}\comp F^j$$ Therefore $\widehat{\hacd}=\{\muj,\;|\;j\in J_{\cd}\}$ are the all linear characters of $\hacd$.
\bt 
Let $\cd\subseteq \cc$ be a fusion subcategory of a pseudo-unitary fusion category $\cc$ with commutative Grothendieck ring. Then we have the following:
\bne
\item
The first orthogonality relation for $\mtc H(\cd)$:
{\Small
\beq\label{hfor}
\sum_{t\in {\cccd}}\frac{\fp(R_t)}{\fp(\xit)^2}\mu_k(\chit)\mu_l(\chitstar)=\delta_{l,k}\frac{\fpcc}{\fp(\cc^k)},\;\text{for any}\; k,l \in {\mtc J_{\cd}},
\eeq
}
where $R_t=\sum_{\ch_i \in \mtc D_t}d_i\ch_i$.
\item 
The second orthogonality relation for $\mtc H (\cd)$:
{\Small
\beq\label{hsor}
\sum_{ k=0}^r \fp(\cc^k) \mu_k({\chit}) \mu_k(\chisstar)=\delta_{s,t}\frac{\fp(\xit)\fp(\xis)\fpcc}{\fp(R_t)},\;\text{for all}\;t, s\in \cccd.
\eeq
}
\ene
\et
\bpf 
Since $\cc$ is a pseudo-unitary fusion category one has $\lam_\cd=\frac{R_\cd}{\fp(R_\cd)}$ for any fusion subcategory and $\dimcc=\fpcc$.  Moreover $\czcc$ is also pseudo-unitary and therefore $\dim(\cc^k)=\fp(\cc^k)$ for any conjugacy class of $\cc$.

Let $j\in \mtc J_\cd$.
If $X_i\in t$ for some $t\in \cccd$  by applying $\muj$ to Equation \eqref{initcoset} one has
\beq\label{incos}
\frac{\muj(\ch_i)}{\fp(X_i)}=\frac{\muj(\chit)}{\fp(\xit)}=\muj(\frac{R_t}{\fp(R_t)}).
\eeq
Note that since $j \in {\mtc J_{\cd}}$ Equation \eqref{lamcsl} implies that $\mu_j(\lamcd)=1$.
 
The second orthogonality relation  for $\cc$ from Equation \eqref{somu} can be written as
$$
\sum_{t \in {\cccd}}\big(\sum_{X_i \in t}\mu_l(\ch_i)\mu_k(\ch_{i^*})\big)=\delta_{l,k}\frac{\fpcc}{\fp(\cc^k)}.
$$
If $l,k\in \mtc J_{\cd}$, using Equation \eqref{incos} one obtains  Equation \eqref{hfor} from above. Indeed,
{\Small
\begin{eqnarray*}
\delta_{l,k}\frac{\fpcc}{\fp(\cc^k)}&=&\sum_{t \in {\cccd}}\sum_{X_i \in t}\mu_l(\ch_i)\mu_k(\ch_{i^*})
\\ &=& \sum_{t \in {\cccd}}\sum_{X_i \in t}\frac{\fp(X_i)\mu_l(\chit)}{\fp(\xit)}\frac{\fp(X_i)\mu_k(\chitstar)}{\fp(\xitstar)}
\\ &=& \sum_{t \in {\cccd}}\frac{\mu_k(\chitstar)\mu_l(\chit)}{\fp(\xit)^2}\big(\sum_{X_i \in t}\fp(X_i)^2\big)=
\\ &=& \sum_{t \in {\cccd}}\frac{\mu_k(\chitstar)\mu_l(\chit)}{\fp(\xit)^2}\fp(R_t)
\end{eqnarray*}
}
The second orthogonality relation follows as usually from the fact that  for complex matrices $(BA=I_n\implies AB=I_n)$. Indeed, denote
$$
a_{tk}:=\sqrt{\frac{\fp(\cc^k)\fp(R_t)}{\fpcc}}\frac{\mu_k(\chit)}{\fp(\xit)},\;k\in \mtc J_\cd, t\in \cccd
$$
and 
$$
b_{lt}:=\sqrt{\frac{\fp(\cc^l)\fp(R_t)}{\fpcc}}\frac{\mu_l(\chitstar)}{\fp(\xitstar)},\;l\in \mtc J_\cd, t\in \cccd.
$$
With these notations it is easy to see that the equality $BA=I_n$ (where $n:=|\cccd|=|\mtc J_\cd|$) is equivalent to the first orthogonality relations. Then the second orthogonality relations are equivalent to $AB=I_n$.

Indeed, 
{\Small
\begin{eqnarray*}
\delta_{s,t}&=&\sum_{k=0}^m a_{tk}b_{ks}\\ &=&\sum_{k=0}^m \sqrt{\frac{\fp(\cc^k)\fp(R_t)}{\fpcc}}\frac{\mu_k(\chit)}{\fp(\xit)}\sqrt{\frac{\fp(\cc^k)\fp(R_s)}{\fpcc}}\frac{\mu_k(\chisstar)}{\fp(\xisstar)}
\\ &=& \frac{\sqrt{\fp(R_t)\fp(R_s)}}{\fp(\xit)\fp(\xis)\fpcc} \sum_{k=0}^m\fp(\cc^k)\mu_k(\chit)\mu_k(\chisstar)
\end{eqnarray*}
}
\epf
\bc\label{div} 
Let $\cd$ be a fusion subcategory of a pseudo-unitary fusion category $\cc$.
With the above notations one has:
\bne
\item For any right coset $t\in \cccd$ and any representative $\xit\in t$ one has 
\beq\label{div1}
\frac{\fp(\xit)^2\fpcc}{\fp(R_t)}\in \mathbb A.
\eeq
\item
Suppose that $\cd\subseteq \ccpt$ acts freely on the simple objects of $\cc$. Then 
\beq\label{div2}\frac{\fpcc}{\fp(\cd) \fp(\cc^j)}\in \mathbb A
\eeq 
for all $j\in \mtc J_{\cd}$.
\ene
\ec
\newcommand{\xim}{X_m}

\bpf
\bne
\item
It follows from Equation \eqref{hsor} for $s=t$ since $\fp(\cc^k)\in \mathbb A$ and $\mu_k(\ch_t)\in \mathbb A$.
\item
Note that in this case, since $\cd$ acts freely, one has $\fp(R_m)=\fp(\cd)\fp(\xim)^2.$ Then equation \eqref{hfor} can be written as
{\Small
\beqn
\sum_{t\in {\cccd}}\mu_k(\chit)\mu_l(\chitstar)=\delta_{l,k}\frac{\fpcc}{\fp(\cd)\fp(\cc^k)},\;\text{for any}\; k,l \in {\mtc J_{\cd}}.
\eeqn
}
Put $k=l$ above, and since $\mu_k(\ch_t)\in \mathbb A$ the divisibility result follows.
\ene
\epf
\bl \label{capcoset}
 Suppose that $\cd,\ca$ are fusion subcategories of a fusion category $\cc$. If $m \in \cccd$ is a right coset of $\cc$ with respect to $\cd$ then the set
 $\ca\cap m$, if not empty, 
is a right coset $\ca$  with respect to $\ca\cap\cd$. 
 \el
 \bpf
Let $X, Y\in \ca \cap m$ be any two simple objects. We have to show that $X\simeq^r_{\ca\cap\cd} Y$ as objects of $\ca$. Since $X, Y$ are in the same right coset with respect to $\cd$, by definition,  there is $Z\in \cd$ such that $m( X, Y\ot Z)>0$. This implies $m(Z, Y^*\ot X)>0$ which in turn gives that $Z\in \ca$ since both $X, Y\in \ca$. Thus $Z\in \cd\cap \ca$. It follows that $X, Y$ are in the same coset of $\ca$ with respect to $\ca\cap \cd$. 
 \epf
\section{Divisibility results for premodular categories}\label{mc}
\subsection{Definition of the braided partition function $M$}
Let $\cc$ be a premodular (i.e. braided and spherical) fusion category. We keep all the notations from the previous section, in particular $\irr(\cc)=\{X_0, X_1,\dots,X_m\}$ and $d_i:=\dim(X_i)$ for all $i$. By \cite[Example 6.14]{scalg}  there is $\comp$-algebra map $\fq: {\cfcc}\ra \cecc$ given by the following formula:
\beq\label{sh}
\fq(\ch_i)=\sum_{i'=0}^m\frac{s_{ii'}}{d_{i'}}E_{i'}.
\eeq
where $S=(s_{ij})$ is the $S$-matrix of $\cc$.

Given a fusion subcategory $\cd$  of a braided fusion category  $\cc$,  the notion of {\it M\"uger centralizer of $\cd$ } was introduced in \cite{dgno2}. Two objects $X$ and $Y$ of $\cc$ centralize each other if $ c_{X, Y}c_{Y, X}=\mtr{id}_{X\ot Y}$.

The centralizer $\cd'$ is defined as the fusion subcategory $\cd'$ of $\cc$ generated by all simple objects $Y$ of $\cc$ centralizing any object objects $X$ of $\cd$.  In the premodular case $X_i$ centralizes $X_{i'}$ if and only if $s_{ii'}=d_id_{i'}$, (see also \cite{proclond}). 
In particular, the centralizer $\cc'$ of $\cc$ is also denoted by $\cz_2(\cc)$ and it is called the M\"uger center of $\cc$. 

Since $\fq(F_j)$ is an idempotent element of $\cecc$ one may write:
\beqn
\fq(F_{j})=\sum_{i \in \mtc A_{j}}E_{i}
\eeqn
for some subset $\ca_j\subseteq \{0, \dots , m\}$. Note that the set $\mtc A_{j}$ might be empty precisely when $\fq(F_{j})=0$. Denote by $\mtc J_2\subseteq \mtc J:=\{0,1,\dots, m\}$ the set of all indices $j$ with $\mtc A_{j}$ not a empty set. Since $\fq(1)=1$ we obtain in this way a partition for the set  of all irreducible representations $\irrcc=\bigsqcup_{j\in \mtc J_2}\mtc {\tilde A}_{j}$ where, to be precise, $\mtc {\tilde A}_j=\{[X_i]\;|\; i\in \mtc A_j\}$.  For any index $0\leq i\leq m$ we denoted by $M(i)$ the unique index $j \in \mtc J_2$ such that $i \in \mtc A_{j}$. One obtains a (unique) function 
 $M:\{0, 1\dots, m\}\ra \mtc J_2$ with the property that $E_i\fq(F_{M(i)})\neq 0$ for all $i \in \{0, 1\dots, m\}$.
%

Since $\cfcc$ is a semisimple commutative algebra it follows that $F_j$ form a $\comp$-linear basis for $\cfcc$. Then for any irreducible character $\ch_i=\sumjtom \alij F_j$ for some $\alij \in\comp$.

The following lemma generalizes \cite[Equation (22)]{CW2}.
\bl
With the above notations, for all $0\leq i,i'\leq m$ one has
\beq\label{symfus}
\frac{\al_{ _{iM(i')}}}{d_i}=\frac{s_{ii'}}{d_id_{i'}}=\frac{\al_{ _{i'M(i)}}}{d_{i'}}.
\eeq
\el
\bpf 
Since $\fq(F_j)=0$ if $j \notin \mtc J$ it follows that 
$
\fq(\ch_i)=\sumjtom \alij \fq(F_j)=\sum_{j \in \mtc J_2}\sum_{i' \in \ca_j}\al_{ _{iM(i')}}E_{i'}=\sum_{i'=0}^r\al_{ _{iM(i')}}E_{i'}.
$
Thus
\beq\label{phirchi}
\fq(\ch_i)=\sum_{i'=0}^m\al_{ _{iM(i')}}E_{i'}.
\eeq
for all $i$. Comparing with Equation \eqref{sh}, \;\text{for all indices } $i$ and $i'$ one has
$
\al_{ _{iM(i')}}=\frac{s_{ii'}}{d_{i'}},$
Now Equation \eqref{symfus} follows since  $s_{ii'}=s_{i'i}$.  
\epf
Recall from \cite[Proposition 8.13.11]{EGNO15}, that for any simple object $X_i\in \irrcc$ there is a character  $\psi_{ _{[X_i]}}:\cfcc\ra \comp$ given by 
$\psi_{ _{[X_i]}}(\ch_{i'})=\frac{s_{ii'}}{d_i}$.

Using  Equation \eqref{symfus} it can be easily seen that
\beq\label{mumi}
\psi_{[X_i]}=\mu_{ _{M(i)}},\text{for all}\;i\in \mtc I.
\eeq 
Indeed, for all indices $i'$ one has:
$$
\psi_{ _{[X_i]}}(\ch_{i'})=\frac{s_{ii'}}{d_{i}}=\al_{ _{i'M(i)}}=\mu_{ _{M(i)}}(\ch_{i'}).
$$
Next Theorem is a generalization of the corresponding result for semisimple Hopf algebras obtained in \cite[Theorem 4.3]{CW2}.
\bt\label{phir1} With the above notations one has 
\beq\label{phirch}
\fq(\ch_i)=\frac{d_i}{\dim(\cc^{M(i)})}\mtr C_{ _{M(i)}} 
\eeq
\et
\bpf
Let $\ch_i=\sum_{j=0}^m\al_{ij}F_j$ as above. Using Equation \eqref{phirchi} and Equation \eqref{symfus} it follows that
\beq\label{intermfus}
\fq(\ch_i)=\sum_{i'=0}^m\al_{ _{iM(i')}}E_{i'}=\sum_{i'=0}^m\frac{\al_{ _{i'M(i)}}d_i}{d_{i'}}E_{i'}=d_i(\sum_{i'=0}^m\frac{\al_{ _{i'M(i)}}}{d_{i'}}E_{i'}).
\eeq
Note that by \cite[Equations (4.3) and (4.8)]{ccc-march}  one has
\beq\label{cjc}
\mtr C_j=
{\dimccj}(\sum_{i'=0}^m \frac{1}{d_{i'} }\al_{ _{i'j}} E_{i'})
\eeq
Then for $j=M(i)$ Equation \eqref{intermfus} becomes $\fq(\ch_i)=\frac{d_i}{\dim(\cc^{M(i)})}\mtr C_{ _{M(i)}} $.
\epf
\subsection{Cosets with respect to $\cz_2(\cc)$}
Recall the notion of cosets with respect to a fusion subcategory from Section \ref{hecke}.

 \bt\label{maincoset}
 Two simple objects $X_i, C_{i'}$ of a pseudo premodular fusion category $\cc$ are in the same coset with respect to $\cz_2(\cc)$ if and only if $M(i)=M(i')$.
 \et
 \bpf
Let $\ch_i:=\ch(X_i)$ and $\ch_{\ip}:=\ch(X_{i'})$. Suppose that the two characters $\ch_i$ and $\ch_{i'}$ are in the same coset with respect to $\cz_2(\cc)$. Then as above from Equation \eqref{cformula} one has:
 $$\frac{\ch_i\lam_{\cc'}}{d_i}=\frac{\ch_{i'}\lam_{\cc'}}{d_{i'}}.$$
 On the other hand $\fq(\lam_{\cc'})=1$ by \cite[Corollary 5.8]{ccc-march}, and applying $\fq$ to the above Equation one has:
 $$
 \fq(\frac{\ch_i\lam_{\cc'}}{d_i})=\frac{1}{\dim(\cc^{M(i)})}\mtr C_{ _{M(i)}}= \fq(\frac{\ch_{i'}\lam_{\cc'}}{d_{i'}})=\frac{1}{\dim(\cc^{M(i')})}\mtr C_{ _{M(i')}}
 $$
which proves that $M(i)=M(i')$.
 
Therefore every set $\mtc {\widetilde A}_j:=\{X_i\;\big|\; M(i)=j\}$ is a union of right cosets with respect to $\cc'$. Clearly, by its definition, the number of the non-empty sets $\ca_j$ equals the cardinality of $\mtc J_2$. 

Note that the relation $\fq(\lam_{\cz_2(\cc}) )=1$ implies also that $\mtc J_2=\mtc J_{\cz_2(\cc)}$. Then from Proposition \eqref{dimhecke} it follows that
 \beq\label{dimz2c}
 \dim \mtc H(\cz_2(\cc))=\big|\mtc J_2\big|.
 \eeq
On the other hand, Equation \eqref{dimz2c} implies that the number of right cosets  also equals the cardinal of $\mtc J_2$. Thus each $\ca_j$ consists of a single right coset with respect to $\cc'$ and the proof is finished.
 \epf
For any $j \in \mtc J_2$ we denote by $\mtc R_j:=\{X_i\;|\; M(i)=j\}$. The above theorem implies that $\mtc R_j$ with $j\in \mtc J_2$ are exactly the right cosets of $\cc$ with respect to $\cz_2(\cc)$.  Denote also ${\mtr R}_j:=\sum\limits_{X_i \in \mtc R_j}d_i\ch_i\in \cfcc$ the regular part of their characters.
\subsection{On the dimension of the cosets}
For any premodular fusion category $\cd$ of $\cc$ we denote by $\mtc { R}(\cd)_j:=\cd\cap \mtc R_j$, the intersection of $\cd$ with each coset $\mtc R_j$ of $\cc$ with respect to $\cc'$. We denote also their regular parts by 
$$
{\mtr R}(\cd)_j:=\sum\limits_{X_i \in \mtc {\mtr R}(\cd)_j}d_i\ch_i\in \cfcc.
$$
\bp\label{proport} Let $\cc$ be a premodular fusion category and $\cd$ a fusion subcategory of $\cc$.
With the above notations it follows that
\beq
 \mtc J_{{\cdp}}=\{M(i)\;|\; \ch_i\in \irr(\cd)\}
\eeq
\text{and}
\beq\label{new}
\dim({\mtr R}(\cd)_j)=\dim(\cd\cap {\ccp})\dim(\cc^j),\; \text{for all} \;j\in \mtc J_\cdp.
\eeq
\ep

\bpf 
Since $\cc$ is spherical one has $d_i=d_{i^*}$ for all $i$. Recall by \cite[Equation (6.8)]{scalg} that 
$
\lam_{\cd}=\frac{1}{\dim(\cd)}\big(\mathlarger{\sum}_{\ch_i\in \irr(\cd)}d_{i^*}\ch_i\big)
$
and Theorem \ref{phirch} gives
$$
f_Q(\lam_{\cd})=\frac{1}{\dim(\cd)}\bigg(\sum_{\ch_i\in \irr(\cd)}\frac{d_id_i^*}{\dim(\cc^{M(i)})}\mtr  C_{ _{M(i)}}\bigg).
$$
Applying now the Fourier transform to the last equality, since $\mtr C_j:=\mtc F^{-1}(F_j)$ it follows that
\begin{eqnarray*}
\mtc F(f_Q(\lam_{\cd}))&=&\frac{1}{\dim(\cd)}\bigg(\sum_{\ch_i\in \irr(\cd)}\frac{d_id_{i^*}}{\dim(\cc^{M(i)})}  F_{ _{M(i)}}\bigg)
\\&=&\frac{1}{\dim(\cd)}\bigg(\mathlarger{\sum}\limits_{j \in \mtc J_2}\big(\mathlarger{\sum}\limits_{\{\ch_i\in \irrcd\;|\; M(i)=j\}}d_id_{i^*} \big)\frac{ F_{j}}{\dim(\cc^{j})}\bigg).
\\ &=& 
\frac{1}{\dim(\cd)}
\mathlarger{\sum}\limits_{\{j \in \mtc J_2 \;|\; j=M(i), \;\ch_i\in \irrcd \}}\frac{\dim({\mtr R}(\cd)_j) }{\dim(\cc^{j})}F_{j}.
\end{eqnarray*}
since, from its definition, one has
$\dim({\mtr R}(\cd)_j):=\bigg(\mathlarger{\sum}\limits_{\{\ch_i\in \irrcd\;|\; M(i)=j\}}d_i^2 \bigg)$.
On the other hand by  \cite[Theorem 1]{ccc-march} one can write
$$
\mtc F(f_Q(\lam_{\cd}))=\frac{\dim({\cdp})}{\dimcc}\lam_{{\cdp}}=\frac{\dim({\cdp})}{\dimcc}\bigg(\mathlarger{\sum}\limits_{j\in \mtc J_{\cd'}} F_j\bigg)
$$
Note that the sphericality of $\cc$ implies by \cite[Prop 4.8.4]{EGNO15} that the coefficients of $F_j$ in both above formulae for $\mtc F(\fq(\lamcd))$ are non-zero scalars. Comparing the coefficients of $F_j$ in the above two formulae for $\mtc F(f_Q(\lam_{\cd}))$ it follows that
$\mtc J_{{\cdp}}=\{M(i)\;|\;\ch_i\in \irr(\cd)\}$ and 
\[
\dim({\mtr R}(\cd)_j)={\dim(\cc^{j})}\frac{\dim({\cdp}){\dim(\cd)}}{\dimcc}
\]
which implies the result since by  \cite[Theorem 3.10]{dgno2} one has
\beq\label{centdim}
\dim(\cd)\dim({\cdp})=\dim(\cc)\dim(\cd\cap {\ccp}).
\eeq
\epf
\bc
With the above notations it follows that
\beq\label{5-mostgen}
\frac{\dimcc \dim(\cc'\cap\cd)}{{\dim({\mtr R}(\cd)_j)}}\in \mathbb A
\eeq
for all $j\in J_{\cd'}$.
\ec
\bpf
Note that equation \eqref{new} gives that $\frac{\dimcc}{\dim(\cc^j)}=\frac{\dimcc \dim(\cc'\cap\cd)}{{\dim({\mtr R}(\cd)_j)}}$.
Since $\ccj$ is a simple object of $\czcc$ the result follows.
\epf


\bc\label{11}
If $\cc$ is an integral braided fusion category of free square dimension and  $\cd\cap \cz_2(\cc)=\vect$ then $\cd$ is pointed. 
\ec

\newcommand{\cdcapccp}{{\cd\cap \cz_2(\cc)}}

\bc \label{disc}
With the above notations it follows that 
\beq\label{newr}
\dim({\mtr R}_j)=\dim({\cc'})\dim(\cc^j),\;\text{for all}\; j \in \mtc J_2.
\eeq
\ec
\bpf
It is Equation \eqref{new} for $\cd=\cc$.
\epf
\subsection{Proof of Theorem \ref{squaredim}}
\bp\label{sumcoset}
Let $\cd$ be a fusion subcategory of a pseudo-unitary premodular category $\cc$. With the above notations one has:
$$
\cd=\bigoplus_{j\in \mtc J_{\cd'}}\mtc {\mtr R}(\cd)_j
$$
is the decomposition of $\cd$ in cosets with respect to $\cd\cap\cc'$.
\ep
\bpf 
By Lemma \ref{capcoset} each non-empty set $\mtc {\mtr R}(\cd)_j$ is a coset of $\cd$ with respect to $\cd\cap \cz_2(\cc)$. On the other hand we will show that
\beqn
\sum\limits_{j\in \mtc J_{\cd'}}\fp(\mtc {\mtr R}(\cd)_j)=\fp(\cd).
\eeqn
This proves that  $\mtc {\mtr R}(\cd)_j$ with $j \in \mtc J_{\cd'}$ are all the right cosets. Indeed, adding all the Equations \eqref{new} for $j\in \mtc J_{\cd'}$ one obtains that
$$
\sum_{j\in \mtc J_{\cd'}}\fp(\mtc {\mtr R}(\cd)_j)=\fp(\cd\cap \cc')\big(\sum
\limits_{j\in \mtc J_{\cd'}}\fp(\cc^j)\big).
$$
On the other hand, since $\cc$ is pseudo-unitary note that by Equation \ref{sumdim} one has 
\beq\label{7-ccc-march}
\sum\limits_{j\in \mtc J_{\cd'}}\fp(\cc^j)=\frac{\fpcc}{\fp(\cd')}=\frac{\fp(\cd)}{\fp(\cd\cap\cz_2(\cc))}
\eeq
and Equation \eqref{centdim} gives that $\sum\limits_{j\in \mtc J_{\cd'}}\fp(\mtc {\mtr R}(\cd)_j)=\fp(\cd).$
\epf

%
 Next one can obtain a  proof of Theorem \ref{squaredim} as follows:
\bpf 
By Proposition \ref{sumcoset} we may suppose that $Y\in \mtc \mtc {\mtr R}(\cd)_j$ for some $j \in \mtc J_{\cd'}$.
If $\cd\cap \cz_2(\cc)=\vect$ then $\mtc {\mtr R}(\cd)_j$ is a coset of $\vect$ inside $\cd$ and therefore it consists of a single element $Y$ of $\cd$. Thus, in this case $\fp({\mtr R}(\cd)_j)=\dim(Y)^2$ and the result follows from Equation \eqref{5-mostgen}.
\epf

\subsection{Proof of Theorem \ref{main3}}
In this subsection we prove our main second result:

\bpf
For the first item suppose that $Y=X_i$ and $M(i)=j$. Since $\cc'$ is pointed all the simple objects in $\mtc R_j$ have the same dimension, namely $\fp(Y)$. If  $r_j:=|\mtc R_j|$ and $G_Y$ is the stabilizer of $Y$ under the action of $\cz_2(\cc)$ it follows that 
$$\fp({\mtr R}_j)=r_j\fp(Y)^2=\frac{\fp(\cc')}{|G_Y|}\fp(Y)^2
$$ 
since $r_j|G_Y|=\fp(\cc')$. Then by Proposition \ref{proport} one has:
\beq\label{ccjf}
\fp(\cc^{j})=\frac{\fp(\mtr R_j)}{\fp(\cc')}=\frac{\fp(Y)^2}{|G_Y|}
\eeq
Note that $\frac{\fp(\cc)}{\fp(\cc^j)}\in \mathbb A$ since $\cc^j$ is a simple object of $\czcc$.  Thus
\beq\label{sty}\frac{\fpcc |G_Y|}{\fp(Y)^2}\in \mathbb A
\eeq 
and the first divisibility result follows.
 
If the action of $\cc'$ is free then by the second item of  Corollary \ref{div}
\beqn\frac{\fpcc}{\fp(\cz_2(\cc)) \fp(\cc^j)}\in \mathbb A
\eeqn
for all $j\in \mtc J_{\cz_2(\cc)}$.
Since in this case $|G_Y|=1$  it follows by  Equation \eqref{ccjf}  that $\fp(\cc^j)=\fp(Y)^2$ and therefore 
$$
\frac{\fpcc}{\fp(\cz_2(\cc))\fp(Y)^2}\in \mathbb A.
$$
\epf
\br
Note that for any simple object $X_i$ of $\cc$  Equation \eqref{div1}  gives that
\beq\label{brdiv1}
\frac{\fp(X_i)^2\fp(\cc)}{\fp(\cc')\fp(\cc^{M(i)})}\in \mathbb A.
\eeq
since  $\fp(R_{M(i)})=\fp(\cz_2(\cc))\fp(\cc^{M(i)})$
\er
\br
As mentioned in the introduction, the divisibility of Equation \eqref{freeptdiv} generalizes \cite[Corollary 3.4.]{yu-sd} in the pseudo-unitary settings.
\er
\bibliographystyle{alpha}
\bibliography{ccts}

\begin{thebibliography}{DGNO10}

\bibitem[BB15]{bbr}
A.~Brugui$\grave{e}$res and S.~Burciu.
\newblock \textnormal{On Normal Tensor Functors and Coset Decompositions for
  Fusion Categories}.
\newblock {\em Appl. Categor. Struct.}, 23:591--608, 2015.

\bibitem[Bur20]{ccc-march}
S.~Burciu.
\newblock \textnormal{Conjugacy classes and centralizers for pivotal fusion
  categories}.
\newblock {\em Monatshefte f\"ur Mathematik}, 193(2):13--46, 2020.

\bibitem[CW10]{CW2}
M.~Cohen and S.~Westreich.
\newblock {\em \textnormal{Higman Ideals and Verlinde type Formulas for Hopf
  algebras}}.
\newblock Trends in Mathematics, Springer Basel AG, 2010.

\bibitem[DGNO10]{dgno2}
V.~Drinfeld, S.~Gelaki, D.~Nikshych, and V.~Ostrik.
\newblock \textnormal{On braided fusion categories, I,}.
\newblock {\em Sel. Math. New Ser.}, 16:1--119, 2010.

\bibitem[EGNO15]{EGNO15}
P.~Etingof, S.~Gelaki, D.~Nikshych, and V.~Ostrik.
\newblock {\em \textnormal{Tensor categories}}, volume 205.
\newblock {\it Mathematical Surveys and Monographs}, American Mathematical
  Society, Providence, RI, 2015.

\bibitem[ENO05]{ENO}
P.~Etingof, D.~Nikshych, and V.~Ostrik.
\newblock \textnormal{On fusion categories}.
\newblock {\em Annals of Mathematics}, 162:581--642, 2005.

\bibitem[M{\"u}g03]{proclond}
M.~M{\"u}ger.
\newblock \textnormal{On the structure of modular categories}.
\newblock {\em Proc. Lond. Math. Soc.}, 87:291--308, 2003.

\bibitem[Shi17]{scalg}
K.~Shimizu.
\newblock \textnormal{The monoidal center and the character algebra}.
\newblock {\em Journal of Pure and Applied Algebra}, 221(9):2338--2371, 2017.

\bibitem[Yu20]{yu-sd}
Z.~Yu.
\newblock \textnormal{On slightly degenerate fusion categories}.
\newblock {\em Journal of Algebra}, 559(1):408--431, 2020.

\end{thebibliography}
\ed

\section*{More to do}
\bne
\item Get the divisibility results from Yu and the orthogonality with unnatural partition-done.
\item
Apply $\fq$ for any coset $t$ of any $\cd$ and see where it goes.
\item 
Use the blms-relations and divisibility 
$$\dimcc c^k_{ij}\in \mathbb A$$
\item Use Harrison for quotient rings.
\ene

\section*{Variuos formulae}

\subsection{$\cd=\cc'$}

Write these general relation in the case $\cd=\cc'$ and then switch between $j's$.

In this case $M(i_t)=t$ gives 
$$\mu_t=\mu_{\xit}$$

 \blue{\Small 
Recall
{\large
\beq\label{symfusp}
\frac{\al_{ _{iM(i')}}}{d_i}=\frac{s_{ii'}}{d_id_{i'}}=\frac{\al_{ _{i'M(i)}}}{d_{i'}}.
\eeq
}
This gives
{\large
\beq\label{symfus2}
\mu_{M(i')}(\ch_i)={\al_{ _{iM(i')}}^m}=\frac{s_{ii'}}{d_{i'}}=\frac{d_i}{d_{i'}}\al_{ _{i'M(i)}}^m=\frac{d_i}{d_{'i}}\mu_{M(i)}(\ch_{i'}).
\eeq
}
}
One has for $i'=i_j$ in Equation 
\green
{\eqref{symfus2}
{\large
$$\mu_{M(i_j)}(\ch_i)=\al_{ _{iM(i_j)}}=\alij=\mu_j(\ch_i)=\frac{d_i}{d_{i_j}}\mu_{M(i)}(\ch_{i_j})=\frac{d_i}{d_{i_j}}\al_{ _{{i_j}M(i)}}=\frac{s_{ii_j}}{d_{i_j}}
$$
}
}
\newpage
\subsection{Image of $\lam_{\cd'}$ by $\fq$}
\blue{\Small 
\br
$$\dim(R_{\cd_j})=\dim(R_\cd(0){\mtr R}(\cd)_j\dim(\cc^j), \; \text{ for all}\; j\in \mtc J_\cdp.
$$
$$\dim(\cd_j)=\dim(\cd_0)\dim(\cc^j), \; \text{ for all}\; j\in \mtc J_\cdp.
$$
\er
}
This shows that if $\ce\subseteq \cd$ then
$$
\frac{\dim(\cd\cap \mtc R_j)}{\dim(\ce\cap \mtc R_j)}=\frac{\dim(\cd\cap \cc')}{\dim(\cce\cap \cc')}
$$
\bp
Let $\cd$ be a fusion subcategory of a premodular fusion category $\cc$. With the above notations one has:
\beq\label{iirev} \irr(\cd')=\bigsqcup_{j\in \mtc J_\cd\cap \mtc J_2}\mtc R_j. 
\eeq
\beq\label{lamcdpr}\fq(\lam_{\cd'})=\sum_{j \in \mtc J_\cd\cap \mtc J_2}\frac{\dim(\mtr R_j)}{\dim(\cd')\dim(\cc^{j})}\mtr C^j=\frac{\dim(\cc')}{\dim(\cd')}\ell_{\cd\vee\cc'}.
\eeq
\ep
By  \cite[Corollary 5.8]{ccc-march}
one has $\fq(\lam_{\cc'})=1_\cecc$. 
\bpf
From the Proof of \cite[Theorem 1]{ccc-march} it follows that
\beq\label{lamcdpr}
\fq(\lamcd)=\sum_{\ch_i\in \irr(\cd')}E_i.
\eeq

On the other hand, since $\lam_\cd=\sum_{j\in \mtc J_{\cd}}F_j$, 
it follows that
$$\fq(\lam_{\cd})=\sum_{j \in J_\cd}\fq(F_j)=\sum_{j \in \mtc J_\cd\cap J_2}\sum_{i \in \ca_j}E_i.
$$
Therefore
$$ \irr(\cd')=\{\ch_i\;|\; i\in\ca_j, \;\text{with}\; j\in \mtc J_\cd\cap \mtc J_2\}=\bigsqcup_{j\in \mtc J_\cd\cap \mtc J_2}\mtc R_j.$$

Note that by definition one has
 $\lam_{\cd'}=\frac{1}{\dim(\cd')}\big(\sum_{\ch_i \in \irr(\cd')} d_i\ch_i.\big)$
 Then
\begin{eqnarray*}
 \fq(\lam_{\cd'})&=&\frac{1}{\dim(\cd')}\bigg(\sum_{\ch_i \in \irr(\cd')} d_{i^*}\fq(\ch_i)\bigg)
 \\&=& \frac{1}{\dim(\cd')}\bigg(\sum_{\ch_i \in \irr(\cd')} d_{i^*}d_i\frac{\mtr C^{M(i)}}{\dim(\cc^{M(i)})}\bigg)
 \\&=& \frac{1}{\dim(\cd')}\bigg(\sum_{j \in \mtc J_\cd\cap \mtc J_2}\big(\sum_{\ch_i \in \mtc R_j}d_{i^*}d_i\big)\frac{\mtr C^{j}}{\dim(\cc^{j})}\bigg)
\\&=&\sum_{j \in \mtc J_\cd\cap \mtc J_2}\frac{\dim(\mtc R_j)}{\dim(\cd')\dim(\cc^{j})}\mtr C^j. 
\end{eqnarray*}
\epf
Using the coset decomposition formula $\dim(R_j)=\dim(\cc')\dim(\cc^j)$ and the formulae becomes:
\beq\label{lamcdpr2}\fq(\lam_{\cd'})=\frac{\dim(\cc')}{\dim(\cd')}\big(\sum_{j \in \mtc J_\cd\cap \mtc J_2}\mtr C^j\big)=\frac{\dim(\cc')}{\dim(\cd')}\ell_{\cd\vee\cc'}.
\eeq
This also follows from Theorem 1 \cite{ccc-march} by replacing $\cd$ with $\cd'$.
\\ In particular for $\cd=\cc$ one has $\mtc J_\cc=\{0\}$ and it follows that
$$
\fq(\lam_{\cc'})=\big(\sum_{j \in \mtc J_\cc\cap \mtc J_2}\mtr C^j\big)=\ell_{\cc}=C^0=1.
$$
\subsection{Image of $\fq(\lamcd)$ and formulae for $\fq(F_0)$}
Note that by definition one has:
$$\fq(F_0)=\fq(\lam_\cc)=\fq(\lam_{(\cc')'})=\sum_{i \in \mtc A_0}E_i=\sum_{\ch_i\in \irr(\cc')}E_i.$$
On the other hand, for $\cd=\cc'$ int the above formula it follows that 
$$\fq(F_0)=\fq(\lam_\cc)=\fq(\lam_{(\cc')'})=\frac{\dim(\cc')}{\dim(\cc)}\ell_{\cc'}=\big(\sum_{j \in \mtc J_2}\mtr C^j\big).$$
\subsubsection*{Other consequences}
\bc 
For any fusion subcategory $\cd$ of a ribbon category one has 
$$\dim(\cd'\cap R_j)=\dim(\cc')\dim(\cc^j)=\dim(R_j), \; \text{ for all}\; j\in \mtc J_\cd \cap J_2.$$
\ec

\bpf 
We apply Proposition \ref{proport} for $\cd'$ instead of $\cc$. Note that $\cd'\supseteq \cc'$ and $\cd'\cap \cc'=\cc'$.
Write the above formula for $\cd'$ instead of $\cd$. Note also that since $\cd''=\cd\vee\cc'$ it follows by Proposition \ref{slcup} that $J_{\cd'}=J_{\cd}\cap J_2$.
\epf
Denote by $\mtc R_j$ the coset of $\cc'$ corresponding to $j\in \mtc J_2$. Thus $\mtc R_j:=\{X_i\;|\; M(i)=j\}$. Using Equation \eqref{mumi} it also follows that $\mtc R_j=\{X_i\;|\;\mu_{X_i}=\mu_j\}$.

For any fusion subcategory $\cd$ of $\cc$ we also denote by $\cd(j)=\cd\cap \mtc R_j$ and let
$${\mtr R}(\cd)_j:=\sum_{X_i \in \cd(j)}d_{i^*}\ch_i.$$

\bc
From the same coset formula, it follows that if $\cd_j=\mtc R_j$ for some $j\in \mtc J_{\cd'}$ then $\cd_j=\mtc R_j$ for all $j \in J_\cdp$. This happens if and only if $\cd\supseteq \cc'$
\ec
\bc
$$
\cd=\oplus_{j \in J_{\cd'}}\cd_j
$$
where $\cd_j:=\cd\cap \mtc R_j$. Moreover, these are the cosets of $\cd\cap \cc'$ inside $\cd$.
\ec
\bpf
Just summing over all $j \in J_\cd'$ it follows one gets $\dim(\cd)$ which implies that these are all the cosets that $\cd$ intersects.

Moreover, it is easy to see that if $X, Y\in \cd_j$ then $X\simeq^{\cd\cap \cc'} Y$ since any $Z$ with  $(X, Z\ot Y)>0$ is necessarily in $\cd$, thus in $\cd\cap \cc'$.
\epf
\bc
$$\cd'=\bigoplus_{j \in J_\cd\cap J_2}\mtc R_j$$
\ec

\bne
\item
The formulae for $\cd\cap \cd'$ cosets inside $\cd$.
\item From the formula, if $\cd$ contains a coset $\mtc R_j$ then $\cd$ contains all the cosets $\mtc R_j$ with $j \in J_\cdp$.
\item
$$\cd'=\oplus_{j\in X} \mtc R_j$$ since $\cd'$ contains $\cc'$.
\ene

\br
\bne
\item
$\{X_i\;|\; M(i)=j\}$ is a coset, therefore this gives after reindexing the cosets (define $R_j$) then
$$\dim(R_j)=\dim(\ccp)\dim(\cc^j).$$
\item
If $\ccp$ is pointed then all in $R_j$ have the same dimension and therefore 
$$\dim(R_j)=r_j\dim(X_j)^2=\dim(\ccp)\dim(\cc^j)$$
which shows that
$$
\dim(X_j)^2\bigg|\dim(\ccp)\dim(\cc).
$$
\item 
If $\ccpt$ acts freely it follows that $r_j=\dim(\ccp)$ and therefore $$\dim(X_j)^2=\dim(\cc^j)\bigg|\frac{\dimcc}{\dim(\cc')}$$
by the second item of Corollary \ref{div}.
\ene
\er

\ed
\section{On the dimension of conjugacy classes and the coset intersection}
Note that for any fusion subcategory one has 
$\cc'\subseteq \cd'$ and therefore ${\mtc J_{\cd'}}\subseteq \mtc J_{\ccp}=\mtc J$.
\vsk One has
$$\mtj(\cd \vee \ce)={\mtc J_\cd}\cap \mtjce. $$ 
\vsk
One has ${\mtc J_\cd}p=\mtjccp$ if and only if $\cdp=\ccp$. Then
$$\cd''=\cc=\cd\vee \cc'$$
\bp\label{proport} Let $\cc$ be a ribbon fusion category and $\cd$ a fusion subcategory of $\cc$.
With the above notations, for any \blue{$j\in \mtj(\cdp)\subseteq \mtjccp=\mtc J$} it follows that:
\beq\label{new}
\sum_{\{i\in \cho_j\;|\ch_i\in \irr(\cd)\}}d_i^2=\dim(\cd\cap {\ccp})\dim(\cc^j)\;\text{and}\; {\mtc J_\cd}p=\{M(i)\;|\; \ch_i\in \irr(\cd)\}.
\eeq
\ep
\bc
This shows that if $\cdp=\ccp$ then $\cd$ intersects all $\cho_j$ with $j \in \mtc J$.
\ec
\vsk One has
$$
{\mtc J_\cd}p=\mtjccp \iff \cdp=\ccp \iff \cd\vee \cc'=\cc\iff
$$
$$
\iff \mtj(\cd \vee \ccp)=\mtj(\cc)=\{0\} \iff \mtj(\cd) \cap \mtj(\ccp)=\{0\}
$$
\bp
One has that $\cdp=\ccp$ fi and only if ${\mtc J_\cd} \cap \mtc J=\{0\}$.
\ep
\vsk
\blue{It follows that
$$
\frac{|\cd_j|}{|\cd_0|}=\frac{|\ce_j|}{|\ce_0|}=\sccj
$$
}
\subsection{When $\cc'$ has a complement?}
$$\cc=\cc' \boxtimes \cd$$
\subsection{Description of the decomposition on the cosets}
Suppose that
$$
\cd=\bigoplus_{j \in \mtj(\ccp)}\cd_j.
$$
Then Proposition 5.22 gives that 
$$
|\cd_j|=\dim(\cd\cap\ccp)\dim(\cc^j)
$$
in particular $\cd_j$ is not empty for any $j\in \mtjccp$.
\subsection*{If $\ccpt\cap \cc'=\vect$}
Then apply coset lemma for $\cd=\ccpt$.

Suppose that $\mtc R_j=\ccpt X_j$. Then $$
\mtc R_j \cap \ccpt=\{\irr (XX_j)\;|\; X\in \ccpt\}
$$
it follows that $\fp(X)=\fp(X_j)$ and $(g, XX_j)=(X^*, X_jg^{-1})=1$, thus $X^*=X_j\ot g^{-1}$ or $X=g\ot X_j^*$, i.e a coset of $\st(X_j)$.
Therefore
$$\dim(\cc^j)=\fp(\mtc R_j\cap \ccpt)=|\st(X_j)|\bigg| \fp(\ccpt).$$
\section{Classification problems}

\section{Relations from M\"uger \cite{proc-lond} and the morphisms $\mu_j$ in the same coset}

\subsection{Morphisms $\muj$ in the ribbon case}
One has that 
$$
\alij=\mu_j(\ch_i),
$$
\blue{\bl
One has $M(\io)=M(\itw)$ if and only if $\mu_{X_\io}=\mu_{X_\itw}$.
\el
\bpf
Use M\"uger's result
\epf
}
\subsubsection{$\muj$ with $j \in \mtc J$ in the ribbon case.}

It follows from [mueger-proc-lond] that
$$
\mu_{X_i}:\cfcc \ra \kk,\;\;Z\sent \frac{s_{ZX_i}}{d(X_i)}
$$ is a morphism of algebras on $\cfcc$. 
\vsk
From the Shimizu-ribbon formula it follows that
$$
\mu_{X_i}=\mu_j, \;\;\text{where} \; j=M(i).
$$
\subsection{The monodromy in the same $\cc'$-coset-rewritten over cosets.}
By Lemma 2.4 i) of [mueger-proc-lond] one has
$$
S(UX,Y)=\frac{1}{d(Y)}S(U,Y)S(X,Y).
$$
If $U\in \irr(\cc')$ then
$$S(UX, Y)=\frac{d(U)d(Y)}{d(Y)}S(X,Y)=d(U)S(X,Y).
$$
In particular for $U=r_{0}$ the normalized regular character of $\cc'$, i.e., the cointegral, one has
$$S(r_0X, Y)=s(X, Y).$$
Since $r_0X=d(X)r_s$ for any $X\in \mtc R_s$ is the coset corresponding to $X$
it follows that
$$S(X,Y)=d(X)S(r_s,Y), \text{for any simple}\;Y
$$
This shows that if $X_{\io}$ and $X_{\itw}$ are in the same coset $\cc'$ then 
\beq\label{eq}
\frac{s_{\io,j}}{d_{\io}}=\frac{s_{i',j}}{d_{\itw}}=\frac{s_{tj}}{d_t}, \;\;\text{for all}\; j\in \mtc J_2(!!)
\eeq
where $X_t$ is a representative  for the coset $\mtc R_t$.
\blue{\bp With the above notations
\beq\label{div}
\frac{d_l^2\dimcc}{\dim(\ccp)\dim(\cc^l)}\in \mathbb A
\eeq
\ep}
Then the second orthogonality relation from Equation \eqref{sos} can be written using the above Equation \eqref{eq} as
\beq\label{sos2}
\sum_{t \in (\cc/\cc')}\sum_{ i\in \mtc R_t}{s_{il}}{s_{i^*l}}=\frac{d_l^2\dimcc}{\dim(\cc^l)}
\eeq
or
\beq\label{sos3}
\sum_{t \in (\cc/\cc')}\frac{s_{tl}s_{t^*l}}{d_t^2}(\sum_{ i\in \mtc R_t}d_i^2)=\frac{d_l^2\dimcc}{\dim(\cc^l)}
\eeq
using the Coset Proposition \ref{proport} and dividing by $\dim(\cc')$ it follows that
\beq\label{sos4}
\sum_{t \in (\cc/\cc')}\frac{s_{tl}s_{t^*l}}{d_t^2}\dim(\mtc \cc^t)=\frac{d_l^2\dimcc}{\dim(\ccp)\dim(\cc^l)}
\eeq
This shows that 
$$
\frac{d_l^2\dimcc}{\dim(\ccp)\dim(\cc^l)}\in \mathbb A
$$
is an algebraic integer.
\subsection{Second orthogonality relation for ribbon category}

See if you get anything new at least in the case of $\ccp \subseteq \ccpt$.

Shimizu relation implies
$$
\fq(\ch_i) =\sum_{\itw=0}^M \frac{s_{i \itw}}{ d_{\itw}} e_{\itw}
$$
On the other hand if $\ch_i=\sumjtom\alij F_j$ it follows 
$$
\fq(\ch_i) = \sum_j \alij (\sum_{\itw \in \ccb_j}e_{\itw})
$$
Thus
\beq
\al_{iM(\itw)} = s_{i \itw}/d_{\itw}
\eeq
For $t,t'$ in the same coset it follows that:
\beq
\al_{iM(t)} =\frac{s_{it}}{d_t}=\al_{iM(t')} =\frac{s_{it'}}{d_{t'}}=\frac{s_{t'i}}{d_{t'}}=\frac{s_{ti}}{d_t}
\eeq
\blue{Divide to the one you change.}

\subsection{Orthogonality relation second}
$$
\sumitom \alij \al_{i^* j} =\frac{\dimcc}{\dimccj}
$$

Suppose that $j=M(t)$. Then the above equation can be written:
$$
\sumitom\al_{iM(t)} \al_{i^* M(t)} =\frac{\dimcc}{\dim(\cc^j)}
$$
or
$$
\sumitom \frac{s_{it}} {d_t}\frac{s_{i^* t}}{d_t}=\frac{\dimcc}{\dim(\cc^{M(t)})}
$$
Recall that
$$
\al_{iM(t)} =\frac{s_{i t}} {d_t} , \;\text{for all}\; M(t)\in \mtc J_2
$$
Written in terms of cosets
$$
\sum_{u\in (\cc/\cc')} \sum_{i \in \co_u} \frac{ d_i s_{ut} }{d_u d_t} \frac{d_{i^*}s_{u^* t}} {d_u d_t} =
$$
Use now 
$$
\sum_{i \in \co_u} d_i^2 = \dim(\cc^u) \dim(\ccp)
$$
\blue{This shows that $\dim\cc'$ divides the norm of each coset.}

It becomes then
$$
(1/d_t^2)\sum_{u\in (\cc/\cc')} \dim(\cc^u) \frac{s_{ut}}{d_u }\frac{s_{u^* t}} {d_u}=\frac{\dimcc}{\dim(\ccp) \dim(\cc^{M(t)})}
$$
This shows that $d_t$ is the same in the range of the coset $M(t)$.

t is fixed and the sum in u is on representatives on cosets.

It is  not missing a coefficient since it does  give the relation from EG88.

Coset relation
$$
s_{it} =\frac{d_i  s_{ut} }{d_u}, \;\text{for all}\; t
$$
if $i$ and $u$ are in the same coset.

Then
$$
\dim(\cc^t)(\dim(\cc')) =r_t d_t^2
$$

\subsection{Relations $\alij$ versus $\sij$}On the other hand, writing as above $\ch_i=\sumjtom \alij F_j$, it follows that 
\beq
f_Q(\ch_i)=\sumjtom \alij f_Q(F_j)=\sum_{j \in \mtc J}\sum_{i' \in \ca_j}\al_{iM(i')}E_i
\eeq
From here it follows that
\beq
\al_{iM(i')}=\frac{s_{ii'}}{d_{i'}}, \;\;\text{from all}\;\;i,i'
\eeq
Since $s_{ii'}=s_{i'i}$ it follows that  
\beq\label{symfus}
\frac{\al_{iM(i')}}{d_i}=\frac{s_{ii'}}{d_id_{i'}}=\frac{\al_{i'M(i)}}{d_{i'}}.
\eeq
generalizing the symmetry Equation (22) from \cite{CW2}.

\subsection{Morphisms of the Grothendieck ring of ribbon fusion categories}
Define $\mu_Y:\cfcc\ra \kk$ 
$$
\mu_Y([X])=\frac{s_{YX}}{d_Y}
$$
This is a morphism of $\gr(\cc)$. 
\vsk
For $Y=X_i$, one has
$$
\mu_{X_i}(X_{\ip})=\frac{s_{i \ip}}{d_i}=\al_{\ip M(i)}.
$$
Therefore $\mu_{X_i}=\mu_{M(i)}$ for all $i$.
\br\label{ch}
One has that $\mu_{X_{\io}}=\mu_{X_{\itw}}$ if and only if they are in the same chunk.
\er
\blue{Understand the categorical proof of the fact it is a morphism from [mueger-proc-lond Lemma 2.3]
}
\subsection{On lemma 2.13 from [mueg-proc]}

Define $p_{\cd}$ as the characteristic function of $\cd$. Then $p_\cd$ is a character of $\gr(\cc)$.
\bl In any fusion category one has
$$
\sum_{j\in \mtj_\cd}\dim(j|\mu_j=d(-)p_\cd[\cc:\cd]
$$
\el

\bpf
One has that $\mu_j(\ch)=\ch(\eta_j)=\frac{\ch(\barcj)}{\sccj}$ and therefore
$$
\sum_{j\in \mtj_\cd}\dim(j|\mu_j(\ch)\numeq{Remark \eqref{ch}}\sum_{j\in \mtj_\cd}\ch(\barcj)\numeq{4.21}\ch(\ell_\cd)\numeq{4.18}d(\ch)p_{\cd}[\cc:\cd]$$
For the last equality we used 4.18 from ccc-march-braided only. For the previous one we used 4.21.
\epf
\subsection{Ribbon case}

$$
\sum_{Y\in \irrcd}d(Y)S(X,Y)=d(X)\dimcd p_{\cdp}(X)
$$
This can be rewritten as
$$<\sum_{Y\in \irrcd}d(Y)^2\mu_Y, X>=\dimcd d(X)p_\cdp(X)
$$ or
$$
\sum_{Y\in \irrcd}d(Y)^2\mu_Y=\dimcd d(-)p_\cdp
$$
Then we can write that
\beq\label{2.13}
\sum_{j \in \mtc J}(\sum_{Y\in \irrcd\cap \ccb^j}d(Y)^2)\mu_j=\dimcd d(-)p_\cdp
\eeq

Using the pervious Proposition it follows that
\beq
\dim(\cd\cap \cc')\left ( \sum_{\{j \in \mtc J, \cd\cap \ccb^j\neq 0\}}\dim(j|\mu_j \right )\numeq{2.13}\dimcd d(-)p_\cdp
\eeq

For $\cd=\cc$ it follows that
\beq\label{cc}
\dim( \cc')\left ( \sum_{\{j \in \mtc J \}}\dim(j|\mu_j \right )=\dimcc d(-) p_\ccp
\eeq
\bc
$Y$ is central if and only if $$\sum_{j\in \mtc J}\dim(j|\muj(Y)=0$$
\ec
For $\cd=\cc'$ it follows that 
\beq\label{ccp}
\dim(\cc')\left ( \sum_{\{j \in \mtc J, \ccpp\cap \ccb^j\neq 0\}}\dim(j|\mu_j \right )=\dimcd d(-)p_\cc
\eeq
\subsection{On Lemma 2.15-[mug-proc]}
\blue{Working with $S$ avoids $f_Q$.}
\bl
$$
\sum_xs(X, Y)s(X, Z)=(\dimcc)\sum_{W\in \ccpp}N^W_{YZ}d_W
$$
\el

\subsection{Slightly non-degenerate fusion categories}
It has a one-dimensional $\psi$ such that $\psi\ch\neq \ch$. Moreover
$\lker_A(\psi)=\phi(A)$.

A coset is $\{\ch, \psi\ch\}$.
\subsection{Standard factorization of rings}
\subsection{On the centrality of chunks}
\blue{Conj1: There are minimal union of chunks that give central characters}
\bpf
??
\epf

The conjecture says that $\mtc J$ can be partitioned in 

\beq
\mtc J=\jc_0\sqcup \jc_1\sqcup\dots \sqcup\jc_r
\eeq
such that 
$$
r_{\jc_s}:=\sum_{\ch_i\in \jc_s}\ch_i(1)\ch_i\in \mtc Z(A^*).
$$
Conjecture 1 allows to introduce the function $M^c$ which is defined by
$$M^c(i)=r,\; \text{if}\; M(i)\in \jc_r.$$
\blue{Conjecture 2:
If an irreducible character $\ch_{\io}$ centralize all the characters in a minimal central part $\jc_l$ then $\jc_l$ is centralized of the entire $\jc_s$ with $s=M^c(\io)$.
}
\noindent
\bne
\item Suppose $\sum_{j \in \mtc C_1}F_j=\sum_{d\in Z_1}\xi_d$ is central in $A^*$. Then 
\beqn
\fq(\sum_{j \in \mtc C_1}F_j)=\sum_{j \in \mtc C_1}\sum_{i \in \ca_j}E_i
\eeqn
is central character of $A^*$.
\item 
\beqn
\fq(F_iA^*)=\fq(F_i)\phi(A)
\eeqn
\blue{It is related with $\phi(A)$-ideals in $A$.}
\ene
\blue{One has that
$$
\oplus_{j \in \mtc C_1}H^*F_j=\oplus\sum_{d\in Z_1} H^*\xi_d
$$
What is $\fq(\xi_d)$?
}
\blue{Never used that $\fq(\ch f)=\fq(\ch)\fq(f)$.}
\vsk
Suppose that $\repapl$ is a normal fusion subcategory. Then $\lam_L$ is a central element of $A^*$ and therefore  $\lam_L$ is sum of minimal central chunks. The question is if the centralizer of such chunk is also central.
\subsection{Generalization of the main result from ccc-march to non factorizable categories}
Where the Fourier transform send various elements.
\subsection{Possible results for the quotients}

\subsection{New files created with new results}
\bne
\item the file of inclusion of algebras ; at the character 
\ene
\subsection{Relation with the quotients}
\bne
\item I have found a formula for $\bar M$. Suppose that
$${\bar F}_{\bar j}=\sum_{j\in \ccb_{\bar j}}F_j$$ 
\noindent 
Then
\blue{
\beqn
\bar \cca_{\bar j}=\cup_{j \in \ccb_{\bar j}} (\cca_j\cap \irr(A//L))
\eeqn
}
\item Work the case when the $A$ is factorizable, i.e. M is identity. In this case $\cca_j=j$. Thus
\blue{
\beqn
\bar \cca_{\bar j}= \ccb_{\bar j}\cap \irr(A//L)
\eeqn
}
\item Furthermore, if $A$ is factorizable then $\rep(A//L)$ is also factorizable if and only if each 
\blue{
$
 \ccb_{\bar j}\cap \irr(A//L)
$
}
has only one element.
\item\red{ In the case both are factorizable the bijection from $A$ restrcits to the bijection from $A//L$. Restricts means that in $\bar \ccb_i$ there is $F_i$ only.
}
\ene
\noindent
Relations of the function $M$ with the quotients. It si related with the other partitions from the other file on inclusion of the character rings.

Then decide when the quotient is factorizable. Let $\pi:A\ra \apl$.

\bl
If $\pi:A\ra B$ is a surjective algebra map of semisimple algebras then
 $\pi(E_i)=\bar {E}_i$ if the module $X_i$ survive the quotient. These covers everything, i.e., 
 $$
 \pi(E_s)=0,\;\text{if}\;\ch_s\notin \irr(B)
 $$
\el
\bpf
Write $${\bar F}_{\bar j}=\sum_{j\in \ccb_{\bar j}}F_j$$
Then
$$
\phi_{\bar R}({\bar F}_{\bar j})=\sum_{\bar i\in {\bar \cca}_{\bar j}}\bar E_i=
$$
$$
=\pi(\sum_{j\in \ccb_{\bar j}}F_j)=\sum_{i \in \cup_{j \in \ccb_{\bar j}} (\cca_j\cap \irr(A//L))}\bar E_i
$$
Therefore
\blue{
\beqn
\bar \cca_{\bar j}=\cup_{j \in \ccb_{\bar j}} (\cca_j\cap \irr(A//L))
\eeqn
}
\epf
\subsection{Conjugacy classes of the quotient} \blue{Those from non-vanishing coefficients (image of $M$) stay simple.}
One has 
$$
\phibr(C_{X_i})={\bar \cc}^{\bar M(i)}
$$
On the other hand from the commutative diagram it follows that 
$$
\phibr(C_{X_i})=\pi(\fq(C_{V_i}))=\pi(\cc^{M(i)})
$$
So for the modules that survives one has 
\blue{
$$
\pi(\cc^{M(i)})={\bar \cc}^{\bar M(i)}
$$
}
\subsection{What happens in general with conjugacy classes}\blue{The other might decompose.}
\vsk
$\pi(\cc^j)$ is also a left coideal and also closed under the left adjoint action of the quotient.
Therefore 
$$
\pi(\cc_j)=\oplus_{\bar j \in \cd_j}\bar \cc^{\bar j}
$$
\red{
This implies that 
$$
\pi(C_j)=\sum_{\bar j\in \ce_j}\eps_j \bar C_{\bar j}
$$
but the sets $\ce_j$ are not disjoint.}
\subsection{Commutativity with Fourier transforms} The diagram applied to $\bar F_{\bar j}$ gives 
that
$$
\bar C_{\bar j}=\frac{\dimkl}{\dimka}(\sum_{j\in \ccb_{\bar j}}\pi(C_j))
$$
\subsection{Other remarks on the idempotents}
The inclusion $(A//L)^*\subset A^*$ and $C(A//L)\subset {\cfcc}$ are both given by $\pi^*$.
\subsection{A remark}
$$
E_i=\frac{d_i}{n}\ch_i(S(\blam_1))\blam_2
$$
After the quotient survives only those characters that are inside $L$, at least if $L$ is a Hopf subalgebra.
\subsection{Few remarks on quasitriangular Hopf algebras}
\bl(L 4.1 Cw-char tabele and n left cid)
$\fq(C)$ is a left coideal stable under the left adjoint action.
\el
\bt
With the above notations one has 
\beqn
\fq(C_{V_i})=\cc^{m(i)} , \text{and},\;\fq(\ch_i)=\al_i \bar C_i
\eeqn
\et
\bpf
Suppose we have proved the second relation $\fq(\ch_i)=\al_i \bar C_i$
\epf
 \section*{Other directions}
 
 \subsection{Connections with $\elcd$}
 Let $\cd$ be a fusion subcategory of a non-degenerate spherical fusion category $\cc$ \blue{with commutative Grotehndieck ring}. We exhibit a central  element $\elcd:=\mtfi(\lam_\cd)\in \cecc$ that completely determines the fusion subcategory $\cd$. 
\mdn
 Suppose that $\cc$ is a non-degenerate  spherical fusion category and $\lam$ is an integral with $\lag\lam, \unu\rag=1$. Then
\beq\label{idmptsums}
\elcd=\frac{\dimcc}{\dimcd}(\sum_{[V_i]\in \irr(\cd)}E_i).
\eeq
\green{
\beq\label{ldldprime-do it}
\ell_{\cd}\ell_{\cd'}=\ell_{\cd\cap \cd'}
\eeq
}
\subsection{Restriction of idempotents from $\czcc$ to $\cc$}
\subsection{The $ \tilde{\mtc F_j}$ for $\czcc$}
\bn{defn}
$\mtc F_j$. $L_j$ is minimal, in the sense that any etale subalgebra of $A$ is  product of this type of algebras.
\end{defn}
\green{\bp
$\ch \in \irr(\cd);\ch F_j=d(\ch)F_j$ then $j\in \mtc J_\cd$.
\ep
\bpf
$\lam_{\cd}=\sum_{j \in \mtc J_\cd}F_j$ and $\lam_\cd\ch=d(\ch)\lam_\cd$. Thus multypling the eigenvector relation with $\lam_\cd$ one obtains
$$\lam_\cd\ch F_j=d(\ch)\lam_\cd F_j$$
which can be rewritten as 
\\
One has the primitive central idempotents in $\cfcd$ $G_i=\sum_{j\in A_i}F_j$
and
$$
\ch=\sum_i\beta_iG_i=\sum_i\beta_i(\sum_jF_j)
$$
The only $i$ for which $\beta_i=d(\ch)$ is $i=0$ corresponding to the integral. 
\epf
}
\bt
On the other hand I know that
$$
\mtc F_j=\rep(\cs_{<<\cc^j>>}).
$$
where $<<\cc^j>>$ is the smallest subalgebra generated by $\cc^j$.
\et
\bpf
Suppose that $\mtc F_j:=\cs_{L_j}$. Since for any $\ch\in \mtc F_j$ one has $\ch f_j=d(\ch)F_j$ it follows that $\cc^j\subseteq L_j$ and therefore $<<\cc^j>>\subseteq L_j$. Therefore $\cs_{<<\cc_j>>}\supseteq \mtc{F}_j$.
\blue{To finish this theorem one needs to use the above Proposition:
}
\epf
\subsubsection{Example-For any factorizable (ribbon) fusion category}
Recall Shimizu's formula which gives that:
$$
\ch_i=\sum_{j}\frac{\sij}{d_j}F_j
$$
Then 
$$
\mtc F_j=\{V_i\;|\; V_i\perp V_j\}=<V_j>'
$$
\subsubsection*{Example -the modular case}

Suppose that
$$
<V_j>=\cs_{\lker_\cc(V_j)}
$$
then 
$$
\mtc F_j=(\cs_{\lker_\cc(V_j)})'=\{\ch_i\;\dim(\cc^i)\subseteq \lker_\cc(V_j)\}.
$$
\subsubsection{Applied for $\czcc$} Denote with $\mtc R$ those $j$ with a nonzero restriction to $\cc$.
For $j \in R$ it follows $V_j=\cc^j$ and that
$$
\tilde{\mtc F}_j=<\cc^j>'
$$
Thus the conjecture $F(\tilde{\mtc F}_j)=\mtc F_j$ is if
$$
F(<\cc^j>)'=\mtc F_j.
$$
On the other hand I know that
$$
\mtc F_j=\rep(H//<<\cc^j>>).
$$

\subsection{Sept 11- further work to do-relations with the quotients}
\bne
\item relations of the braidings with the quotients
\item Is $D(H_8)$ prime?
\ene
\ed